\documentclass[12pt]{amsart}
\usepackage[dvips]{graphicx}
\usepackage{amscd}
\usepackage{amsmath}
\usepackage{amssymb}
\usepackage{latexsym}
\usepackage[english]{babel}

\newtheorem{example}{Example}[section]

\newtheorem{theorem}[example]{Theorem}

\newtheorem{proposition}[example]{Proposition}

\def\Proof{\noindent \it Proof -- \rm}
\def\qed{\hspace{3.5mm} \hfill \vbox{\hrule height 3pt depth 2 pt width 2mm}
\bigskip}

\def\Sym{{\bf Sym}}
\def\PQSym{{\bf PQSym}}        
\def\WQSym{{\bf WQSym}}        
\def\SG{{\mathfrak S}}
\def\SGQSym{{\it \SG QSym}}    
\def\SGSym{{\bf \SG Sym}}      
\def\EFSym{{\bf EFSym}}   
\def\K{{\mathbb K}}
\def\FF{{\mathcal F}}
\def\GG{{\mathcal G}}
\def\PP{{\mathcal S}} 
\def\<{\langle}
\def\>{\rangle}
\def\Std{{\rm std}}
\def\pack{{\rm pack}}
\def\S{{\bf S}}
\def\M{{\bf M}}
\def\card{{\rm card}}

\newcommand{\tun}{\begin{picture}(5,0)(-2,-1)
\put(0,0){\circle*{2}}
\end{picture}}

\newcommand{\tdeux}{\begin{picture}(7,7)(0,-1)
\put(3,0){\circle*{2}}
\put(3,0){\line(0,1){5}}
\put(3,5){\circle*{2}}
\end{picture}}

\newcommand{\ttroisun}{\begin{picture}(15,8)(-5,-1)
\put(3,0){\circle*{2}}
\put(-0.75,0){$\vee$}
\put(6,7){\circle*{2}}
\put(0,7){\circle*{2}}
\end{picture}}
\newcommand{\ttroisdeux}{\begin{picture}(5,12)(-2,-1)
\put(0,0){\circle*{2}}
\put(0,0){\line(0,1){5}}
\put(0,5){\circle*{2}}
\put(0,5){\line(0,1){5}}
\put(0,10){\circle*{2}}
\end{picture}}

\newcommand{\tquatreun}{\begin{picture}(15,12)(-5,-1)
\put(3,0){\circle*{2}}
\put(-0.75,0){$\vee$}
\put(6,7){\circle*{2}}
\put(0,7){\circle*{2}}
\put(3,7){\circle*{2}}
\put(3,0){\line(0,1){7}}
\end{picture}}
\newcommand{\tquatredeux}{\begin{picture}(15,18)(-5,-1)
\put(3,0){\circle*{2}}
\put(-0.75,0){$\vee$}
\put(6,7){\circle*{2}}
\put(0,7){\circle*{2}}
\put(0,14){\circle*{2}}
\put(0,7){\line(0,1){7}}
\end{picture}}
\newcommand{\tquatretrois}{\begin{picture}(15,18)(-5,-1)
\put(3,0){\circle*{2}}
\put(-0.75,0){$\vee$}
\put(6,7){\circle*{2}}
\put(0,7){\circle*{2}}
\put(6,14){\circle*{2}}
\put(6,7){\line(0,1){7}}
\end{picture}}
\newcommand{\tquatrequatre}{\begin{picture}(15,18)(-5,-1)
\put(3,5){\circle*{2}}
\put(-0.75,5){$\vee$}
\put(6,12){\circle*{2}}
\put(0,12){\circle*{2}}
\put(3,0){\circle*{2}}
\put(3,0){\line(0,1){5}}
\end{picture}}
\newcommand{\tquatrecinq}{\begin{picture}(9,19)(-2,-1)
\put(0,0){\circle*{2}}
\put(0,0){\line(0,1){5}}
\put(0,5){\circle*{2}}
\put(0,5){\line(0,1){5}}
\put(0,10){\circle*{2}}
\put(0,10){\line(0,1){5}}
\put(0,15){\circle*{2}}
\end{picture}}

\newcommand{\tdun}[1]{\begin{picture}(10,5)(-2,-1)
\put(0,0){\circle*{2}}
\put(2,-3){\tiny #1}
\end{picture}}

\newcommand{\tddeux}[2]{\begin{picture}(12,5)(-2,-1)
\put(0,0){\circle*{2}}
\put(0,0){\line(0,1){5}}
\put(0,5){\circle*{2}}
\put(2,-3){\tiny #1}
\put(2,3){\tiny #2}
\end{picture}}

\newcommand{\tdtroisun}[3]{\begin{picture}(20,12)(-5,-1)
\put(3,0){\circle*{2}}
\put(-0.75,0){$\vee$}
\put(6,7){\circle*{2}}
\put(0,7){\circle*{2}}
\put(5,-3){\tiny #1}
\put(8,5){\tiny #2}
\put(-6,5){\tiny #3}
\end{picture}}
\newcommand{\tdtroisdeux}[3]{\begin{picture}(12,12)(-2,-1)
\put(0,0){\circle*{2}}
\put(0,0){\line(0,1){5}}
\put(0,5){\circle*{2}}
\put(0,5){\line(0,1){5}}
\put(0,10){\circle*{2}}
\put(2,-3){\tiny #1}
\put(2,3){\tiny #2}
\put(2,9){\tiny #3}
\end{picture}}

\newcommand{\tdquatreun}[4]{\begin{picture}(20,15)(-5,-1)
\put(3,0){\circle*{2}}
\put(-0.6,0){$\vee$}
\put(6,7){\circle*{2}}
\put(0,7){\circle*{2}}
\put(3,7){\circle*{2}}
\put(3,0){\line(0,1){7}}
\put(5,-3){\tiny #1}
\put(8.5,5){\tiny #2}
\put(1,10){\tiny #3}
\put(-5,5){\tiny #4}
\end{picture}}
\newcommand{\tdquatredeux}[4]{\begin{picture}(20,20)(-5,-1)
\put(3,0){\circle*{2}}
\put(-0.75,0){$\vee$}
\put(6,7){\circle*{2}}
\put(0,7){\circle*{2}}
\put(0,14){\circle*{2}}
\put(0,7){\line(0,1){7}}
\put(5,-3){\tiny #1}
\put(8,5){\tiny #2}
\put(-6,5){\tiny #3}
\put(-6,12){\tiny #4}
\end{picture}}
\newcommand{\tdquatretrois}[4]{\begin{picture}(20,20)(-5,-1)
\put(3,0){\circle*{2}}
\put(-0.75,0){$\vee$}
\put(6,7){\circle*{2}}
\put(0,7){\circle*{2}}
\put(6,14){\circle*{2}}
\put(6,7){\line(0,1){7}}
\put(5,-2){\tiny #1}
\put(9,5){\tiny #2}
\put(-5,5){\tiny #4}
\put(9,12){\tiny #3}
\end{picture}}
\newcommand{\tdquatrequatre}[4]{\begin{picture}(20,17)(-5,-1)
\put(3,5){\circle*{2}}
\put(-0.75,5){$\vee$}
\put(6,12){\circle*{2}}
\put(0,12){\circle*{2}}
\put(3,0){\circle*{2}}
\put(3,0){\line(0,1){5}}
\put(5,-3){\tiny #1}
\put(5,3){\tiny #2}
\put(8,12){\tiny #3}
\put(-5,12){\tiny #4}
\end{picture}}

\newcommand{\tdvdeux}[3]{\begin{picture}(15,15)(0,-1)
\put(3,0){\circle*{2}}
\put(1.5,2){.}
\put(1.5,4){.}
\put(1.5,6){.}
\put(3,9){\circle*{2}}
\put(3,9){\line(0,1){5}}
\put(3,14){\circle*{2}}
\put(5,-2){\tiny #1}
\put(6,7){\tiny #2}
\put(6,12){\tiny #3}
\end{picture}}

\newcommand{\tdvquatretrois}[5]{\begin{picture}(20,30)(-5,-1)
\put(3,0){\circle*{2}}
\put(1.5,2){.}
\put(1.5,4){.}
\put(1.5,6){.}
\put(3,9){\circle*{2}}
\put(-0.75,9){$\vee$}
\put(6,16){\circle*{2}}
\put(0,16){\circle*{2}}
\put(6,23){\circle*{2}}
\put(6,16){\line(0,1){7}}
\put(5,-2){\tiny #1}
\put(5,7){\tiny #2}
\put(9,14){\tiny #3}
\put(-5,14){\tiny #5}
\put(9,21){\tiny #4}
\end{picture}}

\input{xy}
\xyoption{all}

\newcommand{\BEQ}{\begin{equation}}     
\newcommand{\BEA}{\begin{eqnarray}}
\newcommand{\EEQ}{\end{equation}}       
\newcommand{\EEA}{\end{eqnarray}}

\newcommand{\Lea}{{\mathrm{Lea}}}
\newcommand{\Roo}{{\mathrm{Roo}}}

\newcommand{\h}{\mathbf{H}}

\newcommand{\T}{\mathbf{T}}
\newcommand{\F}{\mathbf{F}}
\renewcommand{\vec}[1]{\boldsymbol{#1}}
\newcommand{\FQSym}{\mathbf{FQSym}}
\newcommand{\rond}[1]{*++[o][F-]{#1}}

\title{Polynomial realizations of some combinatorial Hopf algebras}

\author[L. Foissy, J.-C.~Novelli, and J.-Y.~Thibon]%
{Lo\"\i c Foissy, Jean-Christophe Novelli and Jean-Yves Thibon}

\address[Foissy]{Laboratoire de Math\'ematiques, Moulin de la Housse \\
Universit\'e de Reims\\ BP 1039, F -- 51687 Reims cedex 2, FRANCE}
\address[Novelli, Thibon]{Institut Gaspard Monge,
Universit\'e Paris-Est Marne-la-Vall\'ee \\
5, Boulevard Descartes \\
Champs-sur-Marne \\ 77454 Marne-la-Vall\'ee cedex 2 \\
FRANCE}
\email[Lo\"\i c Foissy]{loic.foissy@univ-reims.fr}
\email[Jean-Christophe Novelli]{novelli@univ-mlv.fr}
\email[Jean-Yves Thibon]{jyt@univ-mlv.fr}

\date{\today}

\begin{document}

\begin{abstract}
We construct explicit polynomial realizations of some combinatorial Hopf
algebras based on various kind of trees or forests, and some more general
classes of graphs, ranging from the Connes-Kreimer algebra to an algebra of
labelled forests isomorphic to the Hopf algebra of parking functions, and to
a new noncommutative algebra based on endofunctions admitting many interesting
subalgebras and quotients.
\end{abstract}



\maketitle


\section{Introduction}

One knows many examples of Hopf algebras based on various kinds of trees
or forests~\cite{Connes,Fo1,Fo2,LR1,PBT1,PBT2,NTtri}.
Such algebras are increasingly popular, mainly because of their application to
renormalization problems in quantum field theory~\cite{Kreimer,Connes}, but
some of them occured earlier in combinatorics~\cite{Grossman,Grossman2} or in
numerical analysis~\cite{HLW}.

The simplest one, generally known as the Connes-Kreimer algebra~\cite{Connes}, 
is a commutative algebra freely generated by rooted trees, endowed with a
coproduct defined in terms of \emph{admissible cuts}.

This is a basic example of a \emph{combinatorial Hopf algebra}, a heuristic
notion encompassing a large class of graded connected Hopf algebras based on
combinatorial objects, endowed with some extra structure such as distinguished
bases, scalar products or degree-preserving products (called internal
products).
A distinctive feature of combinatorial Hopf algebras is that products and
coproducts in distinguished bases are given by combinatorial algorithms.
However, in many cases, the basis elements can be realized as polynomials
(commutative or not) in some auxiliary set of variables, in such a way that
the product of the algebra becomes the usual product of polynomials, and the
coproduct a simple trick of ``doubling the variables'' (see,
\emph{e.g.},~\cite{NCSF6,NTpark,NTtri} for detailed examples).

Such a construction was not known for the Connes-Kreimer algebra, despite the
fact that it is one of the simplest examples. The present paper provides such
a construction, which in turn will be obtained by specialization of a new
realization of a Hopf algebra of labelled forests~\cite{FU}, itself isomorphic
to the (dual) Hopf algebra of parking functions~\cite{NTpark}. This provides
as well realizations of the noncommutative Connes-Kreimer algebra (isomorphic
to the Loday-Ronco algebra of planar binary trees)~\cite{Fo1,Fo2,LR1}, and new
morphisms between these algebras and other combinatorial Hopf algebras.

Previously known realizations were defined in terms of an auxiliary alphabet
$A$, endowed with some ordering. A given combinatorial Hopf algebra is then
realized by interpreting the elements of some basis as the sum of all words
over $A$ sharing some specific property ({\it e.g.}, descent set,
standardization, packing, parkization), the product is then the ordinary
product of polynomials, and the coproduct is the \emph{ordinal sum} $A+B$
of two isomorphic copies of the ordered set $A$.

As we shall see, it is possible to extend this approach to the algebras of the
Connes-Kreimer family, provided that one replaces the order on $A$ by another
kind of binary relation, for which an analog of the ordinal sum can be
defined.  This construction works for a slightly more general class of graphs,
and we can obtain for example a new Hopf algebra based on endofunctions,
regarded as a generalization of labelled forests where the roots can be
replaced by cycles.

\bigskip
{\bf Acknowledgements.} Partially supported by a PEPS project of the CNRS. 

\section{Rooted trees and rooted forests}

In all the paper, $\K$ will denote a field of characteristic zero.

\subsection{Reminders on rooted trees and forests}

A {\it rooted tree} is a finite tree with a distinguished vertex called the
root.
A rooted forest is a finite graph $\mathcal{F}$ such that any connected
component of $\mathcal{F}$ is a rooted tree.
The set of vertices of the rooted forest $\mathcal{F}$ is denoted by
$V(\mathcal{F})$.

Let $\mathcal{F}$ be a rooted forest. The edges of $\mathcal{F}$ are oriented
downwards (from the leaves to the roots).
If $v,w \in V(\mathcal{F})$, with  $v\neq w$, we shall write $v \rightarrow w$
if there is an edge in $\mathcal{F}$ from $v$ to $w$,
and $v \twoheadrightarrow w$ if there is an oriented path from $v$ to $w$ in
$\mathcal{F}$.

Let $\vec{v}$ be a subset of $V(\mathcal{F})$. We shall say that $\vec{v}$ is
an {\it admissible cut} of $\mathcal{F}$, and we shall write
$\vec{v} \models V(\mathcal{F})$, if $\vec{v}$ is totally disconnected, that
is to say that there is no path from $v$ to $w$ in $\FF$  for any pair
$(v,w)$ of distinct elements of $\vec{v}$.
If $\vec{v} \models V(\mathcal{F})$, we denote by $\Lea_{\vec{v}}\mathcal{F}$
the rooted subforest of $\mathcal{F}$ obtained by keeping only the vertices
above $\vec{v}$, that is to say $\{ w \in V(\mathcal{F}), \: \exists v \in
\vec{v}, \:w \twoheadrightarrow v \}\cup \vec{v}$.  We denote by
$\Roo_{\vec{v}}\mathcal{F}$ the rooted subforest obtained by keeping the other
vertices.

\subsection{The Connes-Kreimer Hopf algebras}

Connes and Kreimer proved in \cite{Connes} that the vector space $\h$
spanned by rooted forests can be turned into a Hopf algebra.
Its product is given by the disjoint union of rooted forests, and the
coproduct is defined for any rooted forest $\mathcal{F}$ by
\begin{equation}
\Delta(\mathcal{F})
= \sum_{\vec{v} \models V(\mathcal{F})}
        \Roo_{\vec{v}}\mathcal{F} \otimes \Lea_{\vec{v}}\mathcal{F}.
\end{equation}
For example,
\begin{equation}
\Delta\left(\tquatredeux\right)=\tquatredeux \otimes 1+1\otimes \tquatredeux
+\ttroisun \otimes \tun+\tdeux \otimes \tdeux+\ttroisdeux \otimes \tun
+\tdeux \otimes \tun\tun+\tun \otimes \tdeux\tun.
\end{equation}
This Hopf algebra is commutative and noncocommutative.
Its dual is the universal enveloping algebra of the free pre-Lie algebra on
one generator~\cite{CL}.

A similar construction can be done on plane forests. The resulting
noncommutative, noncocommutative Hopf algebra $\h_{NCK}$ is called the
noncommutative Connes-Kreimer Hopf algebra~\cite{Fo1,Fo2}. It is isomorphic to
the Hopf algebra of planar binary trees~\cite{LR1}.

\section{Ordered rooted trees and permutations}

We recall here a generalization of the construction of product and the
coproduct of $\h$ to the space spanned by ordered rooted forests introduced
in~\cite{FU}.

\subsection{Hopf algebra of ordered trees}

An \emph{ordered (rooted) forest} is a rooted forest with a total order on the
set of its vertices. The set of ordered forests will be denoted by $\F_o$;
for all $n\geq0$, the set of ordered forests with $n$ vertices will be denoted
by $\F_o(n)$.
An ordered (rooted) tree is a connected ordered forest. The set of ordered
trees will be denoted by $\T_o$; for all $n \geq 1$, the set of ordered trees
with $n$ vertices will be denoted by $\T_o(n)$. The $\K$-vector space
generated by $\F_o$ is denoted by $\h_o$. It is a graded space, the
homogeneous component of degree $n$ being $Vect(\F_o(n))$ for
all $n \in \mathbb{N}$.

For example,
\begin{equation}
\begin{split}
\T_o(1)&=\{\tdun{1}\},\\
\T_o(2)&=\{\tddeux{1}{2},\tddeux{2}{1}\},\\
\T_o(3)&=\left\{\tdtroisun{1}{2}{3},\tdtroisun{2}{1}{3},\tdtroisun{3}{1}{2},
\tdtroisdeux{1}{2}{3},\tdtroisdeux{1}{3}{2},\tdtroisdeux{2}{1}{3},
\tdtroisdeux{2}{3}{1},\tdtroisdeux{3}{1}{2},\tdtroisdeux{3}{2}{1}\right\};
\end{split}
\end{equation}

\begin{equation}
\begin{split}
\F_o(0)&=\{1\},\\
\F_o(1)&=\{\tdun{1}\},\\
\F_o(2)&=\{\tdun{1}\tdun{2},\tddeux{1}{2},\tddeux{2}{1}\},\\
\F_o(3)&=\left\{
  \begin{array}{c}
   \tdun{1}\tdun{2}\tdun{3},
   \tdun{1}\tddeux{2}{3},\tdun{1}\tddeux{3}{2},\tdun{2}\tddeux{1}{3},
   \tdun{2}\tddeux{3}{1},\tdun{3}\tddeux{1}{2},\tdun{3}\tddeux{2}{1}, \\[3mm]
   \tdtroisun{1}{3}{2},\tdtroisun{2}{3}{1},\tdtroisun{3}{2}{1},
   \tdtroisdeux{1}{2}{3},\tdtroisdeux{1}{3}{2},\tdtroisdeux{2}{1}{3},
   \tdtroisdeux{2}{3}{1},\tdtroisdeux{3}{1}{2},\tdtroisdeux{3}{2}{1}
  \end{array}\right\}.
\end{split}
\end{equation}

If $\FF$ and $\GG$ are two ordered forests, then the rooted forest $\FF\GG$ is
seen as the ordered forest such that, for all $v \in V(\FF)$, $w \in V(\GG)$,
$v<w$.
This defines a noncommutative product on the the set of ordered forests. 
For example, the product of $\tdun{1}$ and $\tddeux{1}{2}$ gives
$\tdun{1}\tddeux{2}{3}$, whereas the product of $\tddeux{1}{2}$ and
$\tdun{1}$ gives $\tddeux{1}{2}\tdun{3}=\tdun{3}\tddeux{1}{2}$. This product
is linearly extended to $\h_o$, which in this way becomes a graded algebra.

The number of ordered forests with $n$ vertices is $(n+1)^{n-1}$, which is
also the number of parking functions of length $n$. By definition, $\h_o$ is
free over \emph{irreducible} ordered forests, which are in bijection with
\emph{connected} parking functions. Hence, as an associative algebra $\h_o$ is
isomorphic to the Hopf algebra of parking functions $\PQSym$ introduced
in~\cite{NTpark}.

If $\FF$ is an ordered forest, then any subforest of $\FF$ is also ordered. 
In~\cite{FU}, a coproduct $\Delta:\h_o \longmapsto \h_o\otimes \h_o$ on $\h_o$
has been defined in the following way: for all $\FF \in \F_o$,
\begin{equation}
\Delta(\FF)=\sum_{\vec{v} \models V(\FF)} \Roo_{\vec{v}}\FF\otimes
\Lea_{\vec{v}}\FF.
\end{equation}
As for the Connes-Kreimer Hopf algebra of rooted trees~\cite{Connes}, one can
prove that this coproduct is coassociative, so $\h_o$ is a graded Hopf
algebra. For example,
\begin{equation}
\Delta\left(\tdquatredeux{2}{3}{4}{1}\right)
=\tdquatredeux{2}{3}{4}{1} \otimes 1+1\otimes \tdquatredeux{2}{3}{4}{1}
+\tdtroisun{1}{3}{2} \otimes \tdun{1}+\tddeux{1}{2} \otimes \tddeux{2}{1}
+\tdtroisdeux{2}{3}{1} \otimes \tdun{1} +\tddeux{1}{2} \otimes \tdun{1}\tdun{2}
+\tdun{1} \otimes \tddeux{3}{1}\tdun{2}.
\end{equation}

\begin{theorem}
As a Hopf algebra, $\h_o$ is isomorphic to the graded dual $\PQSym^*$ of
$\PQSym$.
\end{theorem}

{\bf Note.} Actually, $\PQSym$ is self-dual, but as we shall see, $\h_o$
admits $\WQSym$ as a natural quotient rather than as a natural subalgebra,
which is also the case of $\PQSym^*$.

\bigskip
\Proof We shall only give here the main ideas of the proof, see \cite{Foi} for
more details.
Another product, denoted by $\nwarrow$, is defined on the augmentation ideal
of $\h_o$: if $\FF$ and $\GG$ are two ordered forests, $\FF \nwarrow \GG$ is
the ordered forest obtained by grafting $\GG$ shifted by the number of
vertices of $\FF$ on the greatest vertex of $\FF$.
For example,
\begin{equation}
\tddeux{1}{2}\nwarrow\tdun{1}\tdun{2}=\tdquatrequatre{1}{2}{4}{3}
\text{\qquad and\qquad }
\tddeux{2}{1}\nwarrow\tdun{1}\tdun{2}=\tdquatreun{2}{4}{3}{1}.
\end{equation}
This product is associative, and satisfies a certain compatibility with the
product of $\h_o$.
The coproduct of $\h_o$ also splits into two parts, separating the admissible
cuts, according to whether the greatest vertex of $\FF$ is in
$\Roo_{\vec{v}}\FF$ of $\Lea_{\vec{v}}\FF$.
These coproducts make the augmentation ideal of $\h_o$ a dendriform coalgebra,
and there is a certain compatibility (called duplicial) between each product
and each coproduct of $\h_o$, making $\h_o$ what is called in~\cite{Foi} a
$Dup$-$Dend$ bialgebra.
Moreover, the Hopf algebra $\PQSym^*$ is a $Dup$-$Dend$ bialgebra. A rigidity
theorem, similar with the rigidity theorem for bidendriform bialgebra
of~\cite{Foi2}, tells then that a graded, connected $Dup$-$Dend$ bialgebra is
free. As a consequence, as $\h_o$ and $\PQSym\*$ have the same
Poincar\'e-Hilbert series, they are isomorphic as graded
$Dup$-$Dend$ bialgebras, so as graded Hopf algebras.
\qed

\subsection{A realization of $\h_o$}

The Hopf algebra of ordered rooted forests can be realized by explicit
polynomials in an auxiliary alphabet of bi-indexed variables
\begin{equation}
A =\{a_{ij},\ 0\le i<j\}.
\end{equation}
On such an alphabet, we consider the relation
\begin{equation}
a_{ij}\prec a_{kl} \ \Leftrightarrow\ j=k\, .
\end{equation}
We call the pair $(A,\prec)$ a $\prec$-alphabet. This is an analog of
the notion of ordered alphabet used for most other combinatorial Hopf algebras. 
If $(B,\prec)$ is
another $\prec$-alphabet, their $\prec$-sum $A\oplus B$ is defined as their
disjoint union endowed with the $\prec$-relation restricting to the original
ones of $A$ and $B$, and such that $a_{ij}\prec b_{kl}$ for all $i<j$ and
$k<l$.

Let  $\FF$ be an ordered forest with $n$ vertices. We attach to the root of
each tree of $\FF$ a \emph{virtual root}, labeled $0$, $-1$, $-2$ and so on.
There is then a natural bijection from the set of edges of $\FF$ (including
the edges from the roots to the virtual roots) and the vertices of $\FF$,
associating with an edge of $\FF$ its initial vertex.
As the set of vertices of $\FF$ is totally ordered, the set of edges of $\FF$
is then totally ordered.
We shall denote by $e_1<\dots<e_n$ the set of edges of $\FF$.

Let $w=w_1\dots w_n$ be a word of length $n$ over $A$. 
We say that $w$ is $\FF$-\emph{compatible} if, for $k,l \in \{1,\ldots,n\}$
such that the initial vertex of $e_k$ is the terminal vertex of $e_l$,
$w_k \prec w_l$. We then write $w\vdash\FF$. Now define the polynomials
\begin{equation}\label{defSF}
S^\FF(A) = \sum_{w\vdash\FF} w\,.
\end{equation}
For example, let
\begin{equation}\label{exFF}
\FF=\tddeux{3}{2}\tdquatretrois{4}{6}{5}{1}
   =\tdvdeux{-1}{3}{2} \tdvquatretrois{0}{4}{6}{5}{1}.
\end{equation}
Then,
\begin{equation}
\begin{split}
S^{\FF}&=
\sum_{\substack{w_3 \prec w_2 \\w_4 \prec w_1,w_6\\ w_6 \prec w_5}}
w_1w_2w_3w_4w_5w_6\\
&=
\sum_{\substack{i_{-1}<i_3<i_2\\i_0<i_4<i_1\\i_0<i_4<i_6<i_5}}
a_{i_4 i_1}a_{i_3 i_2}a_{i_{-1} i_3}a_{i_0 i_4}a_{i_6 i_5}a_{i_4 i_6}.
\end{split}
\end{equation}

\begin{theorem}
\label{t32}
The polynomials $S^\FF(A)$ provide a realization of $\h_o$. That is,
\begin{equation}
S^\FF S^\GG = S^{\FF\GG}
\end{equation}
and if we allow $A$ and $B$ to commute and identify $P(A)Q(B)$ with
$P\otimes Q$,
\begin{equation}
S^\FF(A\oplus B) = \Delta S^\FF\,.
\end{equation}
\end{theorem}

\Proof
Let us prove that the realization is faithful.
Let $w=a_{i_1,j_1}\ldots a_{i_n,j_n}$ be a word on the alphabet $A$.
Let $J(w)=\{j_1,\ldots,j_n\}$ and $j(w)=\card(J(w))$. Let $i(w)$ be the number
of subscripts $k$ such that $i_k \notin J(w)$.
Then $(i(w),j(w))$ defines a bidegree on $K\<\<A\>\>$. We totally order
$\mathbb{N}^2$ by the lexicographic order.
If $P=\sum x_w w$ is an element of $K\<\<A\>\>$, we denote by $\tilde{P}$
the component of $P$ of maximal $(i(w),j(w))$, if it exists.
In particular, if the degree of the words appearing in $P$ is bounded, then
$\tilde{P}$ exists.
\smallskip

Let $\FF$ be an ordered forest. The degree of the words $w$ appearing in
$S^\FF(A)$ is the number $n$ of vertices of $\FF$, so $\tilde{S}^\FF(A)$
exists. If $w=a_{i_1,j_1}\ldots a_{i_n,j_n}$ appears in $\tilde{S}^\FF$, then
necessarily $j_1,\ldots,j_n$ are all distinct (that is to say $j(w)=n$) and
$i(w)$ is the number of roots of $\FF$.
So $w$ allows to reconstruct $\FF$: there is an edge from the vertex $k$ to
the vertex $l$ in $\FF$ if, and only if, $j_k=i_l$.
As a consequence, the $\tilde{S}^\FF(A)$ are linearly independent.
And so are the $S^\FF(A)$.
\smallskip


Consider now a word $w=c_{i_1j_1}\ldots c_{i_nj_n}$ in $S^\FF(A\oplus B)$. 
If $c_{i_kj_k}$ belongs to $B$, and if $l\rightarrow k$ in $\FF$, then
$c_{i_kj_k} \prec c_{i_lj_l}$, so $c_{i_lj_l}$ also belongs to $B$. As a
consequence, there exists a unique admissible cut $\vec{v}$ such that the
vertices of $\FF$ labelled by those subscripts $k$ such that $c_{i_kj_k}$
belongs to $B$ is $\Lea_{\vec{v}}\FF$ and the vertices of $\FF$ indexed by
those subscripts $k$ such that $c_{i_kj_k}$ belongs to $A$ is
$\Roo_{\vec{v}}\FF$.
Moreover, $w$ is a word appearing in
$S^{\Roo_{\vec{v}}\FF}(A) \otimes S^{\Lea_{\vec{v}}\FF}(B)$.
Conversely, any word appearing in
$S^{\Roo_{\vec{v}}\FF}(A) \otimes S^{\Lea_{\vec{v}}\FF}(B)$ appears in
$S^\FF(A\oplus B)$. So
\begin{equation}
S^\FF(A\oplus B)
=\sum_{\vec{v}\models \FF}
 S^{\Roo_{\vec{v}}\FF}(A) \otimes S^{\Lea_{\vec{v}}\FF}(B)=\Delta S^\FF.
\end{equation}
\qed

For example, let $\FF=\tdquatretrois{$1$}{$3$}{$4$}{$2$}$. Then

\begin{equation}
S^\FF
=\sum_{\substack{w_1\prec w_2,w_3,\\w_3 \prec w_4}}w_1w_2w_3w_4,
\end{equation}
so that
\begin{equation}
S^\FF(A)
=\sum a_{i_0i_1} a_{i_1i_2} a_{i_1i_3} a_{i_3i_4},
\end{equation}
and
\begin{equation}
\begin{split}
S^\FF(A\oplus B)
=& \sum a_{i_0i_1} a_{i_1i_2} a_{i_1i_3} a_{i_3i_4}
  +\sum a_{i_0i_1} b_{j_1j_2} a_{i_1i_3} a_{i_3i_4}
  +\sum a_{i_0i_1} a_{i_1i_2} a_{i_1i_3} b_{j_3j_4}\\
& +\sum a_{i_0i_1} b_{j_1j_2} a_{i_1i_3} b_{j_3j_4}
  +\sum a_{i_0i_1} a_{i_1i_2} b_{j_1j_3} b_{j_3j_4}\\
& +\sum a_{i_0i_1} b_{j_1j_2} b_{j_1j_3} b_{j_3j_4}
  +\sum b_{j_0j_1} b_{j_1j_2} b_{j_1j_3} b_{j_3j_4}\\
=&S^\FF(A)+S^{\tdtroisdeux{1}{2}{3}}(A)S^{\tdun{1}}(B)
  +S^{\tdtroisun{1}{3}{2}}(A)S^{\tdun{1}}(B)\\
& +S^{\tddeux{1}{2}}(A)S^{\tdun{1}\tdun{2}}(B)
  +S^{\tddeux{1}{2}}(A)S^{\tddeux{1}{2}}(B)
  +S^{\tdun{1}}(A)S^{\tdun{1}\tddeux{2}{3}}(B)+S^\FF(B).
\end{split}
\end{equation}

\subsection{Another realization of $\h_o$}
\label{otherrea}

Instead of adding a virtual root to each tree, one can add a loop to the
root. We then have to redefine the alphabet $A$ as
\begin{equation}\label{otherA}
A = \{a_{ij}\ | \ 1\le i\le j \}
\end{equation}
and the $\prec$ relation as
\begin{equation}
a_{ij}\prec a_{jk}\ \text{for $i\le j$ and $j<k$}\,.
\end{equation}
The definition of $\FF$-compatibility remains unchanged: a word $w$ is
$\FF$-compatible if, whenever the initial vertex of edge $e_k$
is the terminal vertex of edge $e_l$, we have $w_k\prec w_l$.
We also keep Definition~(\ref{defSF}) of $S^\FF(A)$. However, note that
the summation indices associated with the virtual roots have disappeared.
For example, the realization of the $S$ indexed by the forest of
Equation~(\ref{exFF}) is
\begin{equation}
S_{\FF} =
\sum_{\substack{i_3<i_2\\i_4<i_1\\i_4<i_6<i_5}}
a_{i_4i_1}a_{i_3i_2}a_{i_3i_3}a_{i_4i_4}a_{i_6i_5}a_{i_4i_6}.
\end{equation}
To recover the coproduct, we redefine the $\prec$-sum of two $\prec$-alphabets
$A$ and $B$ by extending the $\prec$-relations of $A$ and $B$ as follows:
\begin{equation}
a_{jk} \prec b_{ii}\ \text{for all $i,j,k$}.
\end{equation}
Then, we have again the identification
\begin{equation}
\Delta S^\FF \simeq S^\FF(A\oplus B)\,.
\end{equation}
For example, with the forest represented in~(\ref{exFF}), $\Delta S^\FF$ would
contain the sum
\begin{equation}
\sum_{\substack{i_3,i_2\\i_4<i_1\\i_6<i_5}}
a_{i_4i_1}b_{i_2i_2}a_{i_3i_3}a_{i_4i_4}b_{i_6i_5}b_{i_6i_6}.
\end{equation}
corresponding to the cut 
$\tdun{2}\tddeux{6}{5}$.

In the sequel, we shall prefer this version, which has a better
compatibility with known realizations of other combinatorial
Hopf algebras.
We shall also give a third realization of $\h_o$ in Section~\ref{endo-sub}.

\subsection{Epimorphism to $\WQSym$}

Let us recall the definition of $\WQSym$, the Hopf algebra of {\it Word
Quasi-Symmetric functions} (cf. \cite{Hiv,NTtri}).
This algebra has many interpretations, \emph{e.g.}, as the Solomon-Tits
descent algebra~\cite{PS,NTtri}, as a centralizer algebra for a kind of
Schur-Weyl duality~\cite{NPT}, and as an algebra of nonlinear difference
operators~\cite{MNT}.

The \emph{packed word} $u=\pack(w)$ associated with a word $w\in A^*$ (over an
ordered alphabet $A$) is obtained by the following process. If
$b_1<b_2<\ldots <b_r$ are the letters occuring in $w$, $u$ is the image of $w$
by the homomorphism $b_i\mapsto a_i$.
A word $u$ is said to be \emph{packed} if $\pack(u)=u$.
The natural basis of $\WQSym$, which lifts the quasi monomial basis of $QSym$,
is labelled by packed words.
It is defined by
\begin{equation}
\M_u=\sum_{\pack(w)=u}w\,.
\end{equation}
In this basis, the product is given by
\begin{equation}
\M_u\M_v =\sum_{\substack{w=u'v'\\ \pack(u')=u,\,\pack(v')=v}}\M_v\,.
\end{equation}

We shall work with the realization of Section~\ref{otherrea}.
Let $\pi$ be the algebra morphism $a_{ij}\mapsto a_j$ from the free
associative algebra on the $a_{ij}$ of~(\ref{otherA}) to the free
associative algebra over single-indexed letters $a_i$.

\begin{proposition}
$\WQSym$ is a quotient Hopf algebra of $\h_o$:
\begin{equation}
\pi(\h_o)=\WQSym.
\end{equation}
\end{proposition}

\Proof
Let $\FF$ be an ordered forest with $n$ vertices.
A packed word $m=a_1\ldots a_n$ will be said to be $\FF$-admissible if
$i\rightarrow j$ in $\FF$ implies that $a_j < a_i$. Then
\begin{equation}
\pi(\FF)
=\sum_{\mbox{\scriptsize $m$ $\FF$-admissible}}
  \left(\sum_{\pack(w)=m} w\right)
=\sum_{\mbox{\scriptsize $m$ $\FF$-admissible}} \M_m.
\end{equation}
So $\pi(\FF) \in \WQSym$.

Let us prove the surjectivity of $\pi$. We totally order packed words by the
lexicographic order. For any packed word $w=a_1\ldots a_n$, let us construct
an ordered forest $\FF_w$ of degree $n$ such that the smallest packed word
appearing in $\pi(\FF_w)$ is $w$. We proceed by induction on $n$. If $n=1$,
then $w=(1)$ and we take $\FF_{(1)}=\tdun{$1$}$. Let us assume the result for
any packed word with $n-1$ letters. We separate the construction of $\FF_w$
into three cases.
\begin{enumerate}
\item $1=a_1=a_2\leq a_3\ldots \leq a_n$.
We then take $\FF_w=\tdun{$1$} \FF_{a_2\ldots a_n}$.
\item $1=a_1<a_2\leq a_3\leq \ldots \leq a_n$.
Then $a_2=2$. Let $a'_2\ldots a'_n=\pack(a_2\ldots a_n)$, that is
to say $a'_i=a_i-1$ for all $i$.
We then take $\FF_w=B^+(\FF_{a'_2\ldots a'_n})$, that is to say the ordered
tree obtained by adding a root to $\FF_{a'_2\ldots a'_n}$, this root being the
smallest element.
\item The word $a_1\ldots a_n$ is not ordered. There exists $\sigma \in
\SG_n$, such that $a_{\sigma^{-1}(1)}\ldots a_{\sigma^{-1}(n)}$ is ordered.
We then take $\FF_w=\sigma\cdot\FF_{a_{\sigma^{-1}(1)}\ldots
a_{\sigma^{-1}(n)}}$, where $\sigma$ acts by changing the order of the
vertices of $\FF_{a_{\sigma^{-1}(1)}\ldots a_{\sigma^{-1}(n)}}$.
\end{enumerate}
It is not difficult to show that these $\FF_w$ give the result. So
$\pi$ is surjective.
\qed

For example,
\begin{equation}
\begin{split}
\pi(\tdun{1})
 &=\M_{(1)}\\
\pi(\tdun{1}\tdun{2})
 &=\M_{(11)}+\M_{(12)}+\M_{(21)}\\
\pi(\tddeux{1}{2})
 &=\M_{(12)}\\
\pi(\tddeux{2}{1})
 &=\M_{(21)}\\
\pi(\tdtroisun{1}{3}{2})
 &=\M_{(123)}+\M_{(132)}+\M_{(122)}\\
\pi(\tdtroisun{3}{2}{1})
 &=\M_{(231)}+\M_{(321)}+\M_{(221)}.
\end{split}
\end{equation}

\subsection{Embedding of the noncommutative Connes-Kreimer algebra}

Let $\overline{\FF}$ be a plane forest. It can be seen as an
ordered forest, by totally ordering the vertices of $\FF$ "up-left",
that is, by performing a left depth-first traversal of the forest
and numbering each vertex on the first encounter. For example,
\begin{equation}
\tquatretrois\ \tquatredeux \mapsto
\tdquatretrois{1}{3}{4}{2}\quad\tdquatredeux{5}{8}{6}{7}
\end{equation}

\begin{proposition}
The map  $\overline{\FF}\mapsto S^\FF$ is a Hopf embedding of the
noncommutative Connes-Kreimer algebra $\h_{NCK}$ into $\h_o$. 
\end{proposition}

\Proof This is clearly compatible with the product, since shifted concatenation
preserves the planar structure, and with the coproduct which is given on both
sides by admissible cuts, the labeling having been chosen such that in $\h_o$,
the coproduct of the image of a plane forest contains only terms corresponding
to plane forests.
\qed

Thus, we have also a polynomial realization of $\h_{NCK}$. For example,
with the version of Section~\ref{otherrea},
\begin{equation}
\tdeux\tquatretrois \mapsto
S^{\tddeux{1}{2}\tdquatretrois{3}{5}{6}{4}}
=
\sum_{\substack{i_1<i_2\\i_3<i_4\\i_3<i_5<i_6}}
a_{i_1i_1}a_{i_1i_2}a_{i_3i_3}a_{i_3i_4}a_{i_3i_5}a_{i_5i_6}.
\end {equation}

\subsection{Embedding of $\h_{NCK}$ into $\WQSym$}

\begin{theorem}
Let $\pi:\h_o\rightarrow \WQSym$ be the projection induced by
$a_{ij}\mapsto a_j$ in the second realization.
Then the restriction of $\pi$ to $\h_{NCK}$ is injective.
\end{theorem}

\Proof The image of a plane forest can be computed explicitely (see below),
and it then clearly appears that these images are linearly independent.
\qed 

Let $B_+$ denote as usual the operation consisiting in connecting the trees
of a (plane) forest to a common root, and define an endomorphism $b$ of
$\WQSym$ by
\begin{equation}
b(\M_u) = M_{1\cdot u[1]}
\end{equation}
where $\cdot$ is the concatenation and $u[1]$ means shifting by $1$ the
letters of $u$, \emph{e.g.}, $b(\M_{2131})=\M_{13242}$. 

We then have
\begin{equation}
\pi\circ B_+ = b\circ\pi\,.
\end{equation}

An ordered forest $\FF$ can be regarded as a poset $P_\FF$ (the roots being
minimal elements). Identifying a map $f$ from $P_\FF$ 
to some $[m]$ with the word $w_f=f(1)f(2)\cdots f(n)$,
we have
\begin{equation}
\pi(S^\FF) = \sum_{f\in \PP(P_\FF)}\M_{w_f}\,.
\end{equation}
where $\PP(P_\FF)$ is the set of increasing surjections from $P_\FF$ to
some $[m]$.

Now, if we restrict to those $\FF$ which are the canonical labeling of  plane
forests, the lexicographically minimal increasing surjection words $w_f$ are
all distinct and allow to reconstruct $\FF$.

\subsection{The noncommutative Fa\`a di Bruno algebra}

Recall that the Fa\`a di Bruno algebra is the Hopf algebra of
polynomial functions on the group of formal diffeomorphisms of the real line
tangent to the identity~\cite{FG}.
It can be identified with the algebra $Sym$ of symmetric functions, endowed
with the coproduct acting on the complete homeogenous functions $h_n$ by
\begin{equation}
\Delta_1 h_n = \sum_{k=0}^n h_k(X)h_{n-k}((k+1)Y)
\end{equation}
or equivalently,
\begin{equation}
\Delta \sigma_1 =\sum_{n\ge 0} h_n\otimes \sigma_1^{n+1}\,,
\end{equation}
where
\begin{equation}
\sigma_1 := \sum_{n}h_n(X) = \prod_{i\geq1} (1-x_i)^{-1},
\text{\quad and\quad}
\sigma_1(\alpha X) = \sum_n h_n(\alpha X) = \sigma_1(X)^\alpha.
\end{equation}
The noncommutative version \cite{Fo2,BFK} can be identified
with the algebra $\Sym$ of noncommutative symmetric
functions~\cite{NCSF1,NT-Lag}, endowed with the coproduct
\begin{equation}
\Delta_1 S_n = \sum_{k=0}^n S_k(A)S_{n-k}((k+1)B)
\end{equation}
or, again
\begin{equation}
\Delta \sigma_1 =\sum_{n\ge 0} S_n\otimes \sigma_1^{n+1}\,,
\end{equation}
where
\begin{equation}
\sigma_1(A) = \sum_n S_n(A) = \prod_{i\geq1}^{\rightarrow} (1-a_i)^{-1},
\text{\quad and\quad}
\sigma_1(\alpha A) = \sigma(A)^\alpha.
\end{equation}

The Fa\`a di Bruno algebra is known to be a Hopf subalgebra
of the Connes-Kreimer algebra, and in the same way, the noncommutative
version can be embedded in $\h_{NCK}$ \cite{Fo2}. Let
\begin{equation}
U= \sum_{F}F = \frac1{1-V}, \quad V=\sum_{T}T = B_+(U)
\end{equation} 
be the sum of all plane forests and the sum of all plane trees in $\h_{NCK}$.
It is shown in \cite{Fo2} that the square $Z=U^2$ of $U$ has the same coproduct
as $\sigma_1$:
\begin{equation}
\Delta_{NCK} Z =\sum_{n\ge 0} Z_n\otimes Z^{n+1}\,.
\end{equation}
Thus, 
composing the map
$S_n\mapsto Z_n$ with the embedding of $\h_{NCK}$ into $\WQSym$, we obtain an 
embedding of the noncommutative Fa\`a di Bruno algebra.

\subsection{Epimorphism to the original Connes-Kreimer algebra}

If, in the above
realization of $\h_{NCK}$,
 we map $a_{ij}\mapsto x_{ij}$ where the $x_{ij}$ are commuting
indeterminates, we obtain a commutative Hopf algebra which turns out
to be the original Connes-Kreimer algebra.
We can even do this at the level of $\h_o$. With both
realizations,
$S^\FF$ and $S^\GG$ have the same image iff the underlying
unordered forests are the same. Thus the image of $\h_o$ is also the Connes-Kreimer
algebra. 

\begin{proposition}
The map $a_{ij}\mapsto x_{ij}$ provides (for both realizations)
a polynomial realization of the Connes-Kreimer algebra.
\end{proposition}

\Proof The fact that $S^\FF(X)$ depends only on the underlying forests
is clear from the definition. Compatibility with the product and coproduct
is also immediate. The only point which has to be checked is that the
map is surjective. This follows from the same argument as in the proof
of Theorem \ref{t32}.
\qed

The commutative images of the polynomials $S^\FF$ are special cases of
polynomials known in numerical analysis, as well as their coproduct formula
(see, \emph{e.g.}, \cite{HLW}).
More precisely, the specialization of these polynomials $S^\FF$ to the
coefficicients of a finite matrix gives the polynomials associated with each
tree by a Runge-Kutta method (here, with a triangular matrix).
The direct construction of the coproduct in terms of the $\prec$-alphabets
presented here is new.

\subsection{Analog of the Schur basis}

The basis $S^\FF$ is \emph{multiplicative}, in the sense that the product of
two basis elements is again a basis element.
In general, combinatorial Hopf algebras admit several interesting bases, and
such multiplicative bases are generally obtained by summing some other
combinatorial basis along intervals of some order. This is the case here.

There is a natural order on ordered forests with a given number $n$ of
vertices, whose cover relation is $\FF<\FF'$ iff $\FF'$ is obtained from $\FF$
by deleting exactly one edge. In other words, considering the egdes of $\FF$
and $\FF'$ as elements of $\{1,\ldots,n\}^2$, $\FF\leq \FF'$ if, and only if,
the set of edges of $\FF'$ is included in the set of edges of $\FF$.

Let us set
\begin{equation}
R_\FF = \sum_{\GG\le \FF}(-1)^{|E(\FF)|-|E(\GG)|}S^\GG.
\end{equation}
For example,
\begin{equation}
R_{\tddeux{1}{2}\tdun{3}}
=S^{\tddeux{1}{2}\tdun{3}}-S^{\tdtroisun{1}{3}{2}}
-S^{\tdtroisdeux{1}{2}{3}}-S^{\tdtroisdeux{3}{1}{2}}.
\end{equation}

Let $\FF$ be a forest with $k$ vertices and let $I \subseteq \{1,\ldots,k\}$.
The restriction $\FF_{\mid I}$ is the subforest of $\FF$ obtained by taking
all the vertices of $\FF$ which are in $I$ and all the edges between these
vertices. As $I$ is totally ordered (as a part of $\{1,\ldots,k\}$),
$\FF_{\mid I}$ is an ordered forest. Hence

\begin{theorem}
\label{t34}
Let $\FF'$ and $\FF''$ be two ordered forests, with respectively $k'$ and
$k''$ vertices. Then
\begin{equation}
R_{\FF'} R_{\FF''}=\sum R_\FF,
\end{equation}
where the sum is over all ordered forests $\FF$ with $k'+k''$ vertices, such
that $\FF_{\mid \{1,\ldots,k'\}}=\FF'$ and
$\FF_{\mid \{k'+1,\ldots,k'+k''\}}=\FF''$.
\end{theorem}

\Proof Let us define another product $\star$ on $\h_o$, given by the formula
we want to prove.
Let us then compute  $S^{\FF'}\star S^{\FF''}$ for any ordered forests $\FF'$
and $\FF''$. By a M\"obius inversion, for any ordered forest $\GG$,
\begin{equation}
S^{\GG}=\sum_{\GG'\leq \GG} R_{\GG'},
\end{equation}
so that
\begin{equation}
S^{\FF'}\star S^{\FF''}=
\sum_{\GG'\leq \FF',\GG''\leq \FF''} R_{\GG'} \star  R_{\GG''}=\sum R_{\GG},
\end{equation}
where the sum is over all ordered forests $\GG$ with $k'+k''$ vertices, such
that $\GG_{\mid \{1,\ldots,k'\}}\leq \FF'$ and $\GG_{\mid
\{k'+1,\ldots,k'+k''\}}\leq \FF''$.
Such an ordered forest $\GG$ is obtained, first by adding edges between
vertices of $\FF'$ and $\FF''$, then edges between vertices of the two
obtained ordered forests. Equivalently, it can be obtained by adding edges
between vertices of $\FF'\FF''$, so that
\begin{equation}
S^{\FF'}\star S^{\FF''}
=\sum_{\GG \leq \FF'\FF''} R_\GG=S^{\FF'\FF''}=S^{\FF'}S^{\FF''}.
\end{equation}
So, $\star$ is the product of $\h_o$.
\qed

For example,
\begin{equation}
\begin{split}
R_{\tdun{1}\tdun{2}}R_{\tdun{1}}&=
R_{\tdun{1}\tdun{2}\tdun{3}}+R_{\tddeux{1}{2}\tdun{3}}
+R_{\tdun{1}\tddeux{2}{3}}
+R_{\tdun{2}\tddeux{1}{3}}+R_{\tdun{1}\tddeux{3}{2}}
+R_{\tdtroisun{3}{2}{1}}+R_{\tdtroisdeux{1}{3}{2}}
+R_{\tdtroisdeux{2}{3}{1}}.\\
R_{\tdun{1}}R_{\tdun{$1$}\tdun{$2$}}&=
R_{\tdun{1}\tdun{2}\tdun{3}}+R_{\tdun{2}\tddeux{3}{1}}
+R_{\tddeux{2}{1}\tdun{3}}
+R_{\tddeux{1}{2}\tdun{3}}+R_{\tddeux{1}{3}\tdun{2}}
+R_{\tdtroisun{1}{3}{2}}+R_{\tdtroisdeux{2}{1}{3}}
+R_{\tdtroisdeux{3}{1}{2}}.
\end{split}
\end{equation}

Let $\FF$ and $\FF'$ be two ordered forests, equal as rooted forests. There
exists a permutation $\sigma$, such that the ordered forest $\FF^\sigma$
obtained from $\FF$ by permuting the indices by $\sigma$ is equal to $\FF'$.
Then, for any ordered forest $\GG \leq \FF$, as $\GG$ is obtained from $\FF$
by adding some edges, $\GG^\sigma \leq \FF'$.
As a consequence,
the commutative image of $R^{\FF}$ and $R^{\FF'}$ are equal.
For any rooted forest $\overline{\FF}$, we then denote by $R^{\overline{\FF}}$
the image of $R^\FF$ in the Connes-Kreimer algebra, where $\FF$ is any ordered
forest with underlying rooted forest $\overline{\FF}$; this does not depend on
the choice of $\FF$.  These elements form a new basis of the Connes-Kreimer
Hopf algebra.

\medskip
Here are the first examples.
In the Connes-Kreimer Hopf algebra,
\begin{eqnarray*}
R^{\tun}&=&\tun\\
R^{\tdeux}&=&\tdeux\\
R^{\tun\tun}&=&\tun\tun-2\,\tdeux\\
R^{\ttroisdeux}&=&\ttroisdeux\\
R^{\ttroisun}&=&\ttroisun\\
R^{\tun\tdeux}&=&\tun\tdeux-\ttroisun-2\,\ttroisdeux\\
R^{\tun\tun\tun}&=&\tun\tun\tun-6\,\tun\tdeux+3\ttroisun+6\,\ttroisdeux
\end{eqnarray*}
\begin{eqnarray*}
R^{\tquatrecinq}&=&\tquatrecinq\\
R^{\tquatrequatre}&=&\tquatrequatre\\
R^{\tquatredeux}&=&\tquatredeux\\
R^{\tquatreun}&=&\tquatreun\\
R^{\tdeux\tdeux}&=&\tdeux\tdeux-2\tquatredeux-2\,\tquatrecinq\\
R^{\ttroisdeux\tun}&=&\ttroisdeux\tun-\tquatredeux-\tquatrequatre
                      -2\,\tquatrecinq\\
R^{\ttroisun\tun}&=&\ttroisun\tun-\tquatreun-2\tquatredeux-\tquatrequatre\\
R^{\tdeux\tun\tun}&=&\tdeux\tun\tun-2\ttroisun\tun-4\,\ttroisdeux\tun
                      -2\tdeux\tdeux+\tquatreun+6\tquatredeux
                      +3\tquatrequatre+6\,\tquatrecinq\\
R^{\tun\tun\tun\tun}&=&\tun\tun\tun\tun-12\tdeux\tun\tun+12\ttroisun\tun
    +24\,\ttroisdeux\tun+12\tdeux\tdeux-4\tquatreun-24\tquatredeux
    -12\tquatrequatre-24\,\tquatrecinq\,.
\end{eqnarray*}

\section{The cocommutative Hopf algebra on permutations}

Besides the self dual Hopf algebra structure (known as $\FQSym$ or as the
Malvenuto-Reutenauer algebra~\cite{Malvenuto}) on the linear span of all
permutations, there is another one which is cocommutative and noncommutative.
It was first described by Grossman and Larson~\cite{Grossman2} in terms of
heap ordered trees. Several other (non-obviously equivalent) constructions can
be found in~\cite{HNT}. 

The starting point of~\cite{HNT} is a commutative algebra, denoted by
$\SGQSym$, spanned by the polynomials
\begin{equation}
M_\sigma=\sum_{i_1<\ldots<i_n}x_{i_1i_{\sigma(1)}}\cdots x_{i_ni_{\sigma(n)}}
\end{equation}
in commuting indeterminates $x_{ij}$ satisfying $x_{ij}x_{ik}=x_{ik}x_{jk}=0$.
The dual Hopf algebra $\SGSym$ is free over the set of connected permutations,
and the dual basis $\S^\sigma$ of $M_\sigma$ satisfies
\begin{equation}
\label{prodS}
\S^\sigma\S^\tau=\S^{\sigma\bullet\tau}\,,
\end{equation} 
where $\bullet$ denotes shifted concatenation~\cite{HNT}.

\begin{theorem}
\label{t41}
Let $A=\{a_{ij}|i,j\ge 1\}$ endowed with the relation $a_{ij}\prec a_{kl}$
iff $j=k$.
Then, the polynomials 
\begin{equation}
\label{defnewS}
\S^\sigma(A):=\sum_{i_1,\ldots,i_n\ge 1}
a_{i_{\sigma^{-1}(1)}i_1}\cdots
a_{i_{\sigma^{-1}(n)}i_n}
\end{equation}
satisfy (\ref{prodS}), and span a Hopf algebra isomorphic to $\SGSym$ for the
coproduct $\Delta F(A)=F(A\oplus B)$.
\end{theorem}

\Proof
The independence of the $\S^\sigma$ is proved in the same way as for
Theorem~\ref{t32}: indeed, in any $\S^\sigma(A)$ appears a word such that all
subscripts $i_k$ are different and such a word allows one to rebuild $\sigma$.
Moreover, the $\S^\sigma$ defined by (\ref{defnewS}) satisfy
(\ref{prodS}). For the coproduct, observe that $\S^\sigma(A)$ can
alternatively be characterized as the sum of all $\sigma$-compatible words,
defined by the condition:

$w=a_{k_1l_1}\cdots a_{k_nl_n}$ is $\sigma$-compatible iff
\begin{equation}
i=\sigma(j) \Rightarrow a_{k_il_i}\prec a_{k_jl_j}\,.
\end{equation}
Hence, $\S^\sigma(A\oplus B)$ is well-defined, and obtained from
$\S^\sigma(A)$ by splitting the set of cycles of $\sigma$ into two parts in
all possible ways, and replacing $a$'s by $b$'s into one of the parts. This is
exactly the coproduct of the basis $\S^\sigma$ of $\SGSym$ as described
in~\cite{HNT}.
\smallskip

For example, let us consider $\sigma=24513$.
Then
\begin{equation}
\S^\sigma(A)
= \sum a_{i_4i_1}a_{i_1i_2}a_{i_5i_3}a_{i_2i_4}a_{i_3i_5},
\end{equation}
so that
\begin{equation}
\S^\sigma
= \sum_{\substack{w_1\prec w_2 \prec w_4 \prec w_1\\w_3 \prec w_5 \prec w_3}}
 w_1w_2w_3w_4w_5.
\end{equation}
Hence,
\begin{equation}
\begin{split}
\S^\sigma(A\oplus B)
=&\sum a_{i_4i_1}a_{i_1i_2}a_{i_5i_3}a_{i_2i_4}a_{i_3i_5}
 +\sum a_{i_4i_1}a_{i_1i_2}b_{i_5i_3}a_{i_2i_4}b_{i_3i_5}\\
&+\sum b_{i_4i_1}b_{i_1i_2}a_{i_5i_3}b_{i_2i_4}a_{i_3i_5}
 +\sum b_{i_4i_1}b_{i_1i_2}b_{i_5i_3}b_{i_2i_4}b_{i_3i_5}\\
=&\ \S^\sigma(A) + \S^{(231)}(A)\S^{(12)}(B)
  + \S^{(12)}(A)\S^{(231)}(B) + \S^\sigma(B).
\end{split}
\end{equation}

\section{A Hopf algebra of endofunctions}

\subsection{Construction}

The commutative Hopf algebra of permutations of~\cite{HNT} is actually a
subalgebra and a quotient of a commutative algebra based on endofunctions.
There is a similar construction here.

Let $A=\{a_{ij}|i\not = j,\ i,j\ge 1\}$, endowed with the relation
$a_{ij}\prec a_{kl}$ iff $j=k$. 
For a function $f:\ [n]\rightarrow [n]$, let us say that a word
$w=w_1 \cdots w_n$ is $f$-compatible iff
$i\not= j$ and $i=f(j)\Rightarrow  w_i\prec w_j$. Define
\begin{equation}
\S^f(A) :=\sum_{w\ \text{$f$-compatible}}w\,.
\end{equation}
For example, representing a function as the list of its image, if $f=(24352)$,
one has
\begin{equation}
\begin{split}
\S^f =
&\sum_{\substack{w_2 \prec w_1,\\ w_2\prec w_5 \prec w_4 \prec w_2}}
   w_1w_2w_3w_4w_5\\
=&\sum_{i\neq k,j,n, k \neq n,l\neq m}
   a_{ij} a_{ki} a_{lm} a_{nk} a_{in}.
\end{split}
\end{equation}

Note that, as before, these elements are linearly independent: any monomial in
$\S^f$ with as many different subscripts as possible allows one to reconstruct
the relations $w_i<w_j$, and hence the images of $f$ (fixed points being the
missing ones).

\begin{theorem}
The $\S^f$ span a subalgebra of $\K\<A\>$, with
\begin{equation}
\S^f\S^g=\S^{f\bullet g},
\end{equation}
where, again, $\bullet$ denotes the shifted concatenation.
This is a (non-cocommutative) Hopf algebra for the coproduct
$\Delta \S^f=\S^f(A\oplus B)$. 
\end{theorem}

\Proof Similar to the proof of Theorem~\ref{t41}.
\qed

Let us give a description of the coproduct. Let $f:[n] \rightarrow [n]$ and
let $I \subseteq [n]$.
Let $f^I:I\rightarrow I$ be the map satisfying
$f^I(x)=f(x)$ if $f(x) \in I$ and $f^I(x)=x$ otherwise.
If $I$ has cardinality $k$, there exists a unique increasing bijection
$\tau_I:I\rightarrow [k]$;
Then $\Std(f^I):=\tau_I \circ f^I \circ \tau_I^{-1}$.
We shall say that $I$ is an ideal of $f$ and write
$I \models f$ if $f^{-1}(I) \subseteq I$.

One then sees that
\begin{equation}
\Delta(S^f)
=\sum_{I\models f} S^{\Std(f^{[n]\setminus I})} \otimes S^{\Std(f^I)}.
\end{equation}

For example, let us consider $f=(23234)$. Then
\begin{equation}
S^f
=\sum_{\substack{w_2 \prec w_1,w_3\\w_3 \prec w_2,w_4\\w_4\prec w_5}}
  w_1w_2w_3w_4w_5,
\end{equation}
so that
\begin{equation}
S^f(A)=\sum_{\substack{j\not=i,k\\ l\not=k,m}} a_{ji}a_{kj}a_{jk}a_{kl}a_{lm}
\end{equation}
and
\begin{equation}
\begin{split}
S^f(A\oplus B)
=& \sum a_{ji}a_{kj}a_{jk}a_{kl}a_{lm}
  +\sum b_{qp}a_{kj}a_{jk}a_{kl}a_{lm}
  +\sum a_{ji}a_{kj}a_{jk}a_{kl}b_{qp}\\
& +\sum b_{qp}a_{kj}a_{jk}a_{kl}b_{rs}
  +\sum a_{ji}a_{kj}a_{jk}b_{pq}b_{qr}\\
& +\sum b_{qp}a_{kj}a_{jk}b_{rs}b_{st}
  +\sum b_{qp}b_{rq}b_{qr}b_{rs}b_{st}\\
=& S^f(A)+S^{(2123)}(A)S^{(1)}(B)+S^{(2323)}(A)S^{(1)}(B)
  +S^{(212)}(A)S^{(12)}(B)\\
& +S^{(232)}(A)S^{(11)}(B)+S^{(21)}(A)S^{(122)}(B)+S^{(f)}(B).
\end{split}
\end{equation}
Hence
\begin{equation}
\begin{split}
\Delta(S^f)
=& S^f \otimes 1+S^{(2123)}\otimes S^{(1)}+S^{(2323)}\otimes S^{(1)}
  +S^{(212)}\otimes S^{(12)} \\
& +S^{(232)}\otimes S^{(11)} +S^{(21)}\otimes S^{(122)}+1\otimes S^{(f)}.
\end{split}
\end{equation}
Note that the ideals of $f$ are
$\emptyset$, $\{1\}$, $\{5\}$, $\{1,5\}$, $\{4,5\}$, $\{1,4,5\}$,
and $\{1,2,3,4,5\}$.
\smallskip

We shall use a graphical representation of endofunctions. If
$f:[n]\rightarrow [n]$, the vertices of the graph associated with $f$ are the
elements of $[n]$, and there is an edge from $i$ to $j$ iff $f(i)=j$ for all
$i\neq j$. For example, the graph associated with $(23234)$ is
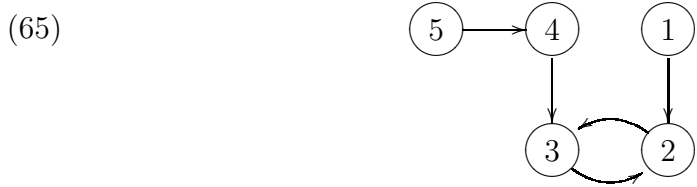
\begin{equation}
\xymatrix{\rond{5}\ar[r]&\rond{4}\ar[d]&\rond{1}\ar[d]\\
&\rond{3}\ar@/_1pc/[r]&\rond{2}\ar@/_1pc/[l]}
\end{equation}
The ideals of $f$ are then given by admissible cuts of the graph (note that
the edges in the cycles cannot be cut).
\smallskip

We shall denote this Hopf algebra of endofunctions by $\EFSym$.

\subsection{Hopf subalgebras and quotients}
\label{endo-sub}

\subsubsection{Permutations}

The $\S^\sigma$, where $\sigma$ runs over permutations, span a Hopf subalgebra
of $\EFSym$ isomorphic to $\SGSym$ (note that the realization is different).
Indeed, if $f$ and $g$ are permutations, then $f \bullet g$ is also a
permutation; if $f$ is a permutation, then its ideals are the disjoint unions
of cycles of $f$, so one recovers the Hopf algebra structure of $\SGSym$
described in  \cite{HNT}.

\subsubsection{Ordered forests}

The $\S^\phi$, where $\phi$ runs over acyclic functions, span a Hopf
subalgebra of $\EFSym$ isomorphic to our first algebra of labelled forests,
hence to $\PQSym^*$. Indeed, if $\FF$ is a labelled forest, we can define an
acyclic function $f_\FF$ in the following way:
if there is an edge from $i$ to $j$ in $\FF$, then $f(j)=i$. If $i$ is a root
of $\FF$, then $f_\FF(i)=i$. For example,
\begin{equation}
\begin{split}
f_{\tdun{1}}&=(1)\\
f_{\tdun{1}\tdun{2}}&=(12)\\
f_{\tddeux{1}{2}}&=(11)\\
f_{\tddeux{2}{1}}&=(22)
\end{split}
\end{equation}
In other words, $f_\FF$ is the endofunction whose graph is $\FF$, the
orientation being implicitly from top to bottom.
Now, $f_{\FF\GG}=f_\FF \bullet f_\GG$. Moreover, the ideals of $f_\FF$ are
the set of the indices $I$ such that the vertices of $\FF$ indexed by $I$ are
a $\Lea_{\vec{v}}\FF$, where $\vec{v}$ runs over the set of admissible cuts of
$\FF$.
So
\begin{equation}
\begin{array}{l}
\h_o \to \EFSym \\
S^\FF \mapsto S^{f_\FF}
\end{array}
\end{equation}
is an injective map of graded Hopf algebras.

\subsubsection{Plane forests}

We have seen that the noncommutative Connes-Kreimer algebra is a Hopf
subalgebra of $\h_o$.
Moreover, for any ordered forest $\FF$, the acyclic function $f_\FF$ is a
nondecreasing parking function if, and only if, $\FF$ is a plane forest.
So the restriction of the embedding $S^\FF \mapsto S^{f_\FF}$ is an
isomorphism from the noncommutative Connes-Kreimer Hopf algebra to the
subspace spanned by the $\S^\pi$, where $\pi$ runs over nondecreasing
parking functions, which is then a Hopf subalgebra of $\EFSym$.
So this gives (yet) another realization of the noncommutative Connes-Kreimer
algebra.

Let $I_o$ be the subspace generated by the $S^f$, where $f$ runs over the set
of endofunctions which are not acyclic. It is clear that $I_o$ is an ideal of
the Hopf algebra of endofunctions. Moreover, if $f$ is not acyclic, then it
contains a cycle $C=i_1\mapsto i_2 \mapsto\ldots \mapsto i_k\mapsto i_1$ of
length $\geq 2$. Let $I$ be an ideal of $f$. If $C \cap I \neq \emptyset$,
then by definition of the ideals, $C \subseteq I$, so $\Std(f^I)$ or
$\Std(f^{[n]-I})$ is not acyclic: this implies that $I_o$ is a Hopf ideal of
the Hopf algebra of endofunctions. So the quotient $\h_o/I_o$ is
isomorphic to the Hopf subalgebra of acyclic endofunctions, so to $\h_o$ and
$\PQSym^*$.

\subsubsection{Nondecreasing functions}

The restriction to nondecreasing functions also gives rise to a Hopf algebra:
if $f$ and $g$ are nondecreasing functions, then $f \bullet g$ is also
nondecreasing, and for any ideal $I$ of $f$, $\Std(f^I)$ and $\Std(f^{[n]-I})$
are also nondecreasing.

\subsubsection{Burnside classes}

The restriction to idempotent functions, or more generally to Burnside classes
($f^p=f^q$) gives, as in the commutative case, Hopf subalgebras: if $f^p=f^q$
and $g^p=g^q$, then $(f\bullet g)^p=(f\bullet g)^q$ and for any part $I$ of
the domain of $f$, $\Std(f^I)^p=\Std(f^I)^q$. Graphically, this corresponds to
endofunctions $f$ such that the graph of $f$ contains only cycles of length
dividing $|p-q|$ and trees of height smaller than $|p-q|$. In particular, for
the idempotent functions, this gives endofunctions whose graph is a corolla,
that is to say a tree of height smaller than~$1$.

\subsubsection{Connes-Kreimer}

The commutative images of the $\S^\phi$ ($a_{ij}\mapsto x_{ij}$) span a
commutative Hopf algebra containing the algebra $Sym$ of ordinary symmetric
functions (as the image of the subalgebra isomorphic to $\SGSym$) and the
Connes-Kreimer Hopf algebra of trees (as the image of the subalgebra
isomorphic to $\h_o$).

\subsection{Analog of the Schur basis}

We define a partial order on the set of endofuctions of $[n]$ for a fixed $n$,
whose cover relation is $f <g$ if there exists an element $j$ of $[n]$
with $f(j)\neq j$, such that $g(k)=f(k)$ if $k\neq j$ and $g(j)=j$. 
For example, for $n=2$, the Hasse graph of this partial order is:

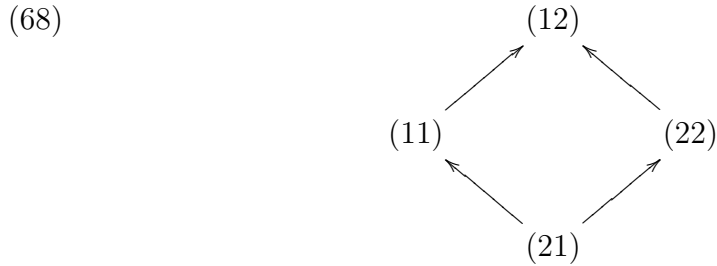
\begin{equation}
\xymatrix{&(12)&\\
(11)\ar[ru]&&(22)\ar[lu]\\
&(21)\ar[ru]\ar[lu]&}
\end{equation}

{\bf Note.} Let $\FF$ and $\GG$ be two ordered forests. It is not difficult to
show that $f_\FF \leq f_\GG$ if, and only if, $\FF \leq \GG$.
\smallskip

For any endofunction $f$, let us set
\begin{equation}
R_f=\sum_{g\leq f} (-1)^{{\rm Fix}(f) - {\rm Fix}(g)} S^g,
\end{equation}
where $Fix(f)$ denotes the number of fixed points of $f$.
By a M\"obius inversion, for any endofunction $f$,
\begin{equation}
S^f=\sum_{g\leq f} R_g.
\end{equation}
By analogy with Theorem~\ref{t34}, one can show

\begin{theorem}
Let $f'$ and $f''$ be two endofunctions of respectively $[n']$ and $[n'']$.
Then
\begin{equation}
R_{f'} R_{f''}=\sum R_f,
\end{equation}
where the sum is over all endofunctions $f$ of $[n'+n'']$,
such that $\Std(f^{[n']})=f'$ and $\Std(f^{[n'+n'']\setminus [n']})=f''$.
\end{theorem}

For example,
\begin{equation}
\begin{split}
R_{(12)}R_{(1)}
=& R_{(121)}+R_{(122)}+R_{(123)}+R_{(131)}+R_{(132)}+R_{(133)}\\
& +R_{(321)}+R_{(322)}+R_{(323)}+R_{(331)}+R_{(332)}+R_{(333)}\\
\end{split}
\end{equation}
\begin{equation}
\begin{split}
R_{(1)}R_{(12)}
=&R_{(111)}+R_{(113)}+R_{(121)}+R_{(123)}+R_{(211)}+R_{(213)}\\
& +R_{(221)}+R_{(223)}+R_{(311)}+R_{(313)}+R_{(321)}+R_{(323)}
\end{split}
\end{equation}
Indeed, for $R_{(12)}R_{(1)}$, one gets all functions such that $f(1)$ is
either $1$ or $3$, and $f(2)$ is either $2$ or $3$, the value $f(3)$ having no
constraint at all.

Let us consider the subspace $I'_o$ generated by the $R_f$, where $f$ runs
over the set of endofunctions $f$ which are not acyclic.
If $f'$ or $f''$ is not acyclic and if $R_f$ appears in $R_{f'}R_{f''}$, then
$f$ is not acyclic.
So $I'_0$ is an ideal. Let us denote by $\overline{R_f}$ the class of $R_f$
modulo $I'_0$.
Note that $\overline{R_f}$ is nonzero if, and only if, $f$ is acyclic,
that is to say there exists an ordered forest $\FF$ such that $f=f_\FF$.
Moreover, if $\FF$ is an ordered forest with $n$ vertices and
$I \subseteq [n]$, then $\FF_{\mid I}=\GG$ if, and only if,
$\Std(f_\FF^I)=f_\GG$. So the map $R_\FF \mapsto \overline R_{f_\FF}$ is an
algebra isomorphism from the algebra $\h_o$ to the algebra of endofunctions
quotiented by $I'_o$.
\smallskip

{\bf Remark.} The ideals $I_o$ and $I'_o$ are different: in degree $2$, $I_o$
is spanned by $S^{(21)}$, whereas $I'_0$ is generated by
$R_{(21)} = S^{(21)}-S^{(11)}-S^{(22)}+S^{(12)}$.



\footnotesize

\begin{thebibliography}{99}
\bibitem{BFK}{\sc C. Brouder}, {\sc A. Frabetti} and {\sc C. Krattenthaler},
{\it Non-commutative Hopf algebra of formal diffeomorphisms}, 
Adv. Math. {\bf 200} (2006), no. 2, 479–524. 
%
\bibitem{CL}{\sc F. Chapoton} and {\sc M. Livernet},
{\it Pre-Lie algebra and the Rooted Trees Operad}, 
Int. Math. Res. Not. IMRN  2001,  no. 8, 395--408.
%
\bibitem{Connes}{\sc A. Connes} and {\sc D. Kreimer},
{\it Hopf algebras, Renormalization and Noncommutative geometry},
Comm. Math. Phys. {\bf 199} (1998), no. 1, 203--242, arXiv:hep-th/9808042.
%
\bibitem{NCSF6}{\sc G. Duchamp, F. Hivert}, and {\sc J.-Y. Thibon},
{\it Noncommutative symmetric functions VI: free quasi-symmetric functions and
related algebras},
Internat. J. Alg. Comput. {\bf 12} (2002), 671--717.
%
\bibitem{FG}{\sc H. Figueroa} and {\sc J.M. Gracia-Bondia},
{\it Combinatorial Hopf algebras in quantum field theory I},
Rev. Math. Phys. {\bf 17} (2005), 881--976.
%
\bibitem{Fo1}{\sc L. Foissy},
{\it Les alg\`ebres de Hopf des arbres enracin\'es d\'ecor\'es} I, 
Bull. Sci. Math.  {\bf 126}  (2002),  no. 3, 193--239. 
%
\bibitem{Fo2}{\sc L. Foissy},
{\it Les alg\`ebres de Hopf des arbres enracin\'es d\'ecor\'es} II, 
Bull. Sci. Math.  {\bf 126}  (2002),  no. 4, 249--288.
%
\bibitem{Foi}{\sc L. Foissy}, 
{\it Ordered forests and parking functions},
Int. Math. Res. Notices (2011) doi: 10.1093/imrn/rnr061,  arXiv:1007.1547.
%
\bibitem{Foi2}{\sc L. Foissy},
{\it Bidendriform bialgebras, trees, and free quasi-symmetric functions},
J. Pure Appl. Algebra  {\bf 209} (2007), no. 2, 439--459,
arXiv:math/0505207.
%
\bibitem{FU}{\sc L. Foissy} and {\sc J. Unterberger},
{\it Ordered forests, permutations and iterated integrals}, preprint,
arXiv:1004.5208.
%
\bibitem{NCSF1}{\sc I.~M. Gelfand, D. Krob, A. Lascoux, B. Leclerc,
V.~S. Retakh} and {\sc J.-Y. Thibon}.
{\it Noncommutative symmetric functions},
Adv. Math. {\bf 112} (1995), 218--348.
%
\bibitem{Grossman}{\sc R. L. Grossman} and {\sc R. G. Larson},
{\it Hopf-algebraic structure of combinatorial objects and differential
operators},
Israel J. Math. {\bf 72} (1990), no. 1--2, 109--117, arXiv:0711.3877.
%
\bibitem{Grossman2}{\sc R. L. Grossman} and {\sc R. G. Larson},
{\it Hopf algebras of heap ordered trees and permutations},
Comm. Algebra {\bf 37} (2009), no. 2, 453--459, arXiv:0706.1327.
%
\bibitem{HLW}{\sc E. Hairer, C. Lubich} and {\sc G. Wanner}, 
{\it Geometric Numerical Integration - Structure-Preserving Algorithms for
Ordinary Differential Equations}, 
Second edition, Springer Series in Computational Mathematics 31.
Springer, Berlin (2006).
%
\bibitem{Hiv}{\sc F. Hivert},
{\it Combinatoire des fonctions quasi-sym\'etriques},
Th\`ese de Doctorat, Marne-La-Vall\'ee, 1999.
%
\bibitem{HNT}{\sc F. Hivert, J.-C. Novelli}, and {\sc J.-Y. Thibon},
{\it Commutative combinatorial Hopf algebras},
J. Algebraic Combin. {\bf 28} (2008), no. 1, 65--95.
%
\bibitem{PBT1} {\sc F. Hivert, J.-C. Novelli}, and {\sc J.-Y. Thibon},
{\it Un analogue du mono\"\i de plaxique pour les arbres binaires de
recherche},
C. R. Acad. Sci. Paris S\'er I Math. {\bf 332} (2002), 577--580.
%
\bibitem{PBT2} {\sc F. Hivert, J.-C. Novelli}, and {\sc J.-Y. Thibon},
{\it The algebra of binary search trees},
Theo. Comp. Sci. {\bf 339} (1) (2005), 129--165.
%
\bibitem{Kreimer}{\sc D. Kreimer}, 
{\it Chen's iterated integral represents the operator product expansion},
Adv. Theor. Math. Phys. {\bf 3} (1999), no. 3, 627--670.
%
\bibitem{LR1}{\sc J.-L. Loday} and {\sc M.O. Ronco},
Hopf algebra of the planar binary trees,
\emph{Adv. Math.} 139 (1998) n. 2, 293--309.
%
\bibitem{Malvenuto}{\sc C. Malvenuto} and {\sc C. Reutenauer},
{\it Duality between quasi-symmetric functions and the Solomon descent
algebra},
J. Algebra {\bf 177} (1995), no. 3, 967--982.
%
\bibitem{MNT}{\sc F. Menous, J.-C. Novelli}, and {\sc J.-Y. Thibon},
{\it Mould calculus, polyhedral cones, and characters of combinatorial Hopf
algebras},
math.CO/11091.1634.
%
\bibitem{NPT}{\sc J.-C. Novelli, F. Patras}, and {\sc J.-Y. Thibon},
{\it Natural endomorphisms of quasi-shuffle Hopf algebras},
math.CO/1101.0725.
%
\bibitem{NTtri}{\sc J.-C. Novelli} and {\sc J.-Y. Thibon},
{\it Polynomial realizations of some trialgebras},
Proc. FPSAC'06 (San Diego, 2006). 
%
\bibitem{NTpark}{\sc J.-C. Novelli} and {\sc J.-Y. Thibon},
{\it Hopf algebras and dendriform structures arising from parking functions},  
Fund. Math. {\bf 193}  (2007),  no. 3, 189--241. 
%
\bibitem{NT-Lag}{\sc J.-C. Novelli} and {\sc J.-Y. Thibon},
{\it Noncommutative Symmetric Functions and Lagrange Inversion}
Adv. App. Math., {\bf 40} (2008), 8--35.
%
\bibitem{NTT}{\sc J.-C. Novelli, J.-Y. Thibon}, and {\sc N.M. Thi\'ery},
{\it Alg\`ebres de Hopf de graphes},
C. R. Acad. Sci., Paris, S\'er. A, {\bf 339} (2004), vol. 9, 607--610.
%
\bibitem{PS}{\sc F. Patras} and {\sc M. Schocker},
{\it Twisted Descent Algebras and the Solomon-Tits Algebra},
Adv. in Math., {\bf 199} (1) (2006), 151--184.
\end{thebibliography}
\end{document}